\input amstex
\documentstyle{amsppt}
\magnification\magstep1         
\baselineskip=10pt
\parindent=0.2truein
\pagewidth{6.6truein} \pageheight{9.0truein} \topmatter \subjclass{Primary:37F45,
Secondary:37F30.}
\topmatter
\endsubjclass
\title Non-persistently recurrent points, qc-surgery and instability of rational
maps with totally disconnected Julia sets
\endtitle
\abstract Let $ R $ be a rational map with totally disconnected Julia set $ J(R). $
If the postcritical set on $ J(R) $ contains a non-persistently recurrent (or
conical) point, then we show that the map $ R $ can not be a structurally stable
map.
\endabstract
\author Peter M. Makienko\endauthor
\address{Permanent address:
Instituto de Matematicas,Av. de Universidad s/N., Col. Lomas de Chamilpa, C.P.
62210,  Cuernavaca, Morelos, Mexico}
\endaddress
\endtopmatter
\rightheadtext{Non-persistently recurrent points, qc-surgery and}
\document
\heading Introduction and Statements
\endheading
Fatou's problem of the density of hyperbolic maps in the space of rational maps is
one of the principal problems in the field of holomorphic dynamics. Due to Man{\'e},
Sad and Sullivan \cite{MSS} we can reformulate this problem in the following way:

{\it If the Julia set $ J(R) $ contains a critical point, then the rational map $ R
$ is a structurally unstable map.}

For convenience we give the definition of the structural stability of a rational
map. For other basic notations and definitions  refer to the book of Milnor
\cite{M}.

\proclaim{Definition 1} Let $ Rat_d $ be the space of all rational maps of degree $
d$ with the topology of  coefficient convergence. A map $ R \in Rat_d $ is called
{\rm structurally stable} if there exists a neighborhood $ U \subset Rat_d $ of $ R
$ such that:

For any map $ R_1 \in U $ there exists a quasiconformal map $ f:\overline{\Bbb
C}\mapsto\overline{\Bbb C} $ conjugating $ R $ to $ R_1$.
\endproclaim
We give a condition "Assumption G" (see below) on the rational map with totally
disconnected Julia set and with a critical point on $ J(R) $ to be unstable. In a
pioneer paper \cite{BH}, Branner and Hubbard prove that the Lebesgue measure of the
Julia set is zero if there exists only one critical point on $ J(R).$ Our result
(see theorem A below) restricted on the Branner-Hubbard case is weaker, but it can
be applied for maps with two or more critical points on $J(R).$
 Let $ R $ be a rational map with totally
disconnected Julia set. Let us normalize $ R $ so that the point $ z = \infty $
becomes the attractive fixed point. Let $ Pc(R) $ be postcritical set of the map $ R
$ and $ P(R) = Pc(R)\cap J(R) $ be postcritical set on the Julia set. Let $ S =
{\Bbb C}\backslash \overline{\cup_nR^{-n}(Pc(R))}, $ then $ R: S\rightarrow S $ is
an unbranched autocovering.

\proclaim{Definition 2}We call a closed simple geodesic $ \gamma \subset S $ {\rm
linked with} $ P(R) $ if the interior $ I(\gamma) $ of $ \gamma $ intersects $
P(R).$
\endproclaim
\proclaim{Assumptions "G"}Let $ R $ be a rational map with totally disconnected
Julia set. Assume there exists a simple closed geodesic $ \gamma \in S $ such that:
\roster \item there exists an infinite subsequence of simple closed geodesics $
\gamma_i \in \cup_n R^{-n}(\gamma) $ linked with postcritical set $P(R),$ and \item
for all $ i =1, ...$ the lengthes $ L(\gamma_i) $ are bounded uniformly away from $
 \infty, $ in the hyperbolic metric on $ S.$
\endroster
\endproclaim

The aim of this paper is to  prove  the following theorem. \proclaim{Theorem A} Let
$ R $ be a rational map with totally disconnected Julia set satisfying Assumption
"G". Then the map $ R $ is not a structurally stable map (that is to say, is an
unstable map.)
\endproclaim
Apriori it is not  clear when the Assumption G holds. We give a series of sufficient
conditions on $ R $ that imply Assumption "G". The next proposition translates
Assumption "G" into the language of "non-persistently recurrent points" on P(R).
\proclaim{Definition 3} A point $ x \in P(R) $ is called {\rm persistently
recurrent} if any backward orbit $ U_0, U_{-1}, ... $ of any neighborhood $ U_0 $ of
$ x $ along $ P(R) $  hits a critical point infinitely many times.
\endproclaim
\proclaim{Lemma (Sufficient condition)} Let $ R $ be a rational map with totally
disconnected Julia set and assume there exists a non-persistently recurrent point $
x \in P(R).$ Then  $ R $ satisfies Assumption "G".
\endproclaim
\demo{Prof}Follows immediately from Definition 3.
\enddemo

Another sufficient condition is connected with the conical points of $ P(R).$
\proclaim{Definition 4} let $ R $ be a rational map, and denote by $ U(x_0,
R^k,\delta) $ the component of $ R^{-k}({\Bbb D}(R^k(x_0),\delta)) $ that contains
$x_0,$ where $ {\Bbb D}(R^k(x_0),\delta) $ is the disk centered in $ R^k(x_0) $ with
radius $ \delta.$ A point $ x_0 $ is {\rm conical} if and only if there is a
constant $ \delta > 0, d \in {\Bbb N}, $ and a sequence $ k_j\to \infty $ such that
$$ R^{k_j}:U(x_0,
R^k,\delta)\mapsto {\Bbb D}(R^k(x_0),\delta)
$$
has degree no more than $ d.$
\endproclaim
Several other notions of conical point appear in the literature. One can see that
the definition of conical point is somehow in the spirit of the notion of conical
set of Lyubich and Minsky \cite{LM}. Definition 4 above appears in \cite{P}, where
Przytycki compares different notions of conical points. McMullen \cite{MM} and
independently Urbanski \cite{DMNU} call a point conical if the mappings in
Definition 4  can be chosen to be conformal.

\proclaim{Theorem 1} Let $ R $ be a rational map with totally disconnected Julia set
and assume there exists a conical point $ x  \in P(R). $ Then $ R $ is an unstable
map.
\endproclaim
The following two results are immediate corollaries of the theorem 1.

\proclaim{Corollary 1} Let $ R $ be a rational map with totally disconnected Julia
set and assume that the postcritical set $ P(R) $ contains a periodic point $ x. $
Then $ R $ is unstable map.
\endproclaim
\demo{Proof} By assumption, the periodic point $ x \in J(R). $ Hence $ x $ is either
parabolic or repelling. Now assume that $ R $ is a structurally stable map, then $ x
$ should be repelling and hence conical. Applying Theorem 1 we are done.
\enddemo
\proclaim{Corollary 2} Let $ R $ be a rational map with totally disconnected Julia
set. Assume $ J(R) = P(R), $ then R is an unstable map.
\endproclaim
\demo{Proof} In this case $ P(R) $ contains all repelling periodic points and by
Corollary 1 we are done.
\enddemo
 \heading Proof of Theorem A
 \endheading

To   prove Theorem A, we use a kind of quasiconformal surgery in the spirit of
Shishikura \cite{Sh}.

Let $ \Delta(r)  $ be a disk of the radius $ r $ centered at $ z = 0 $ and $ \Delta
= \Delta(1). $ Let $ A(p,q) = \{z: p <\vert z\vert < q $ be a ring. Let $ A(p) =
A(p, 1), p < 1 $ and $ a = \frac{1 + 3p}{4} \in A(p, \frac{1 + p}{2}) \subset A(p) $
be a point. Now we define a quasiconformal homeomorphism $ f_p: \Delta \mapsto
\Delta $ as follows: \roster \item $ f_p(z) = \frac{z + a}{1 + az} $ on $ \Delta(p)
$ and \item $ f_p $ is a quasiconformal mollifier on $ A(p), $ that is

i)$ f_p = id $ on $\partial \Delta $ and

ii)$ f_p = \frac{z + a}{1 + az} $ on other boundary component of $ A(p) $ and

iii) the $ L_\infty-$ norm of the dilatation
$\mu=\frac{\overline{\partial}f_p}{\partial f_p} $ is minimal among all dilatations
of the quasiconformal homeomorphisms satisfying  i)-ii).
\endroster
\proclaim{Remark 1} Note that the $ L_\infty-$ norm of  $\mu $ depends only on the
modulus of the ring $ A(p), $ or in  other words, if $ p_i \in \Delta $ converges to
$ p_0 \in \Delta, $ then the $ L_\infty-$ norm of the dilatations $ \mu_i $ are
uniformly bounded away from $ 1. $
\endproclaim

According to results of D. Sullivan (\cite{S}) and C. McMullen, D. Sullivan
(\cite{MS}) the space of full orbits of the points on $ S $ forms a Riemann surface
$ S(R) $ which is (conformally) the torus with a finite number of punctures: the
punctures correspond to the full orbits of the critical points belonging to $ F(R).
$ Hence  there exists a fundamental domain $ F\subset F(R) $ for the action $ R $ on
$ F(R). $ We can choose the fundamental domain as follows: \roster \item $ F \subset
F(R) $ is a closed topological ring, \item the boundary components $ \alpha_1\cup
\alpha_2 = \partial F \subset S $ are smooth closed Jordan curves \item $
R(\alpha_1) =\alpha_2$ and restrictions $ R^n_{\vert F} $ are univalent for all $ n.
$
\endroster

Let  $ O(F) $ be the full orbit of the fundamental domain $ F. $ Let $ \alpha
\subset S $ be any geodesic, then $ \alpha $ intersects a finite number, say $
n(\alpha), $ of elements of $ O(F), $ say $ F_1(\alpha), ..., F_{n(\alpha)}(\alpha).
$

\proclaim{Remark 2} By the properties of the fundamental domain we can always assume
that there exists $ i_0 $ so that the forward orbit $ O_{+}(F_{i_0}(\alpha)) =
\cup_{i\geq 1}R^{i}(F_{i_0}(\alpha)) $ never intersects the interior of the geodesic
$ \alpha. $ For convenience we redefine $ F_1(\alpha) = F_{i_0}(\alpha). $
\endproclaim

Let $ B(\alpha) \subset S\cap \{\cup_{i = 1}^{n(\alpha)} F_i(\alpha) $ be an annulus
containing $ \alpha $ as a non-trivial curve with modulus $ m(\alpha) $ of $
B(\alpha) $ as large as possible. Note that $ B(\alpha) $ is not unique.
 Now let $ \beta \subset  S $ be an iterated preimage of $
\alpha $ (that is there exists an integer $ k $ such that $ R^k(\beta) = \alpha $).
If $ d(\beta) $ is the degree of the covering $ R^k:\beta\mapsto\alpha, $ then the
hyperbolic length $ l(\beta) = d(\beta)l(\alpha), $ and $ m(\beta) \geq
\frac{m(\alpha)}{d(\beta)}, $ as well as $ n(\beta) \leq d(\beta)n(\alpha). $

Let us start with any closed simple geodesic $ \gamma \in S $ linked with $ P(R).$
Now we associate a qc-homeomorphism $ f(\gamma, p): \overline{\Bbb
C}\mapsto\overline{\Bbb C} $ as follows:

Let $ I(\gamma) $ be the interior of $ \gamma $ and the point $ b $ be the first hit
in $ I(\gamma) $ of the forward orbit of a critical point $ c \notin I(\gamma. $
Let $ h : I(\gamma)\mapsto \Delta $ be the Riemann map with $ h(b) = 0, h^\prime(b)
= 1. $ Now let $ p > 0 $ be an number so that $ A(p) \subset h(B(\gamma)). $
Adjusting $ h $ by a rotation we can construct a conformal map $ \phi(\gamma):
I(\gamma)\mapsto \Delta $ so that  the point $ a = \frac{1 + 3p}{4} \in
\phi(F_1(\gamma)). $ Then we set
$$
f(\gamma, p) = \cases \phi(\gamma)^{-1}\circ f_p \circ \phi(\gamma) \text{ on }
I(\gamma),\\
id \text{ off of } I(\gamma)
\endcases
$$

Hence for any simple closed geodesic $ \gamma \subset S $ and a suitable number $ 0
< p < 1 $  we can define a quasi-regular map $ P(\gamma, p) = f(\gamma, p)\circ
R:\overline{\Bbb C}\mapsto\overline{\Bbb C}. $

\proclaim{Lemma 1} Let $ \gamma \subset S $ be closed simple geodesic linked with $
P(R) $ and $ 1 > p > 0 $ be a suitable number. Then \roster\item there exists an
invariant conformal structure $ \sigma $ on $ \overline{\Bbb C} $ so that $
P(\gamma, p): (\overline{\Bbb C}, \sigma)\mapsto(\overline{\Bbb C}, \sigma $ is a
holomorphic map, \item the norm of the dilatation of $ \sigma $ (that is
$L_\infty-$norm of the corresponding Beltrami differential) depends  only on the
numbers $ p $ and $ n(\gamma). $
\endroster
\endproclaim
\demo{Proof}  Follows immediately from the definition of $ f(\gamma, p). $
\enddemo
By the Riemann Mapping Theorem there exists a quasi-conformal homeomorphism   $
f_\sigma $ fixing the points $ 0, 1 $ and $ \infty $  so that $ R(\gamma, p) =
f_\sigma\circ R\circ f_\sigma^{-1} $ is a rational map.

Let $ s(R) $ be the number of critical points whose forward orbits converge to $
\infty. $ \proclaim{Corollary 1} Let $ \gamma  $ and $ p $ be as in Lemma 1 above.
Assume that $ \gamma $ has a sufficiently small spherical diameter. Then $
s(R(\gamma, p) \geq s(R) + 1. $
\endproclaim
\demo{Proof} Let the spherical diameter of $ \gamma $ be so small that the interior
$ I(\gamma) $ does not contain any critical point of the Fatou set $ F(R). $ Hence
if the critical point $ c \in F(R), $ then $ P^n(\gamma, p)(c) = R^n(c) \to\infty. $
Now let $ c \in J(R) $ be the critical point coming from the definition of $
f(\gamma, p). $ Then again by the construction, $ P^n(\gamma, p)(c) \to\infty. $
\enddemo

Now we are ready to prove  Theorem A. Let $ \gamma_i \in \cup_n R^{-n}(\gamma) $
 be coming geodesics from the assumption. Let us redefine: \roster \item $ \phi_i = \phi(\gamma_i), $
 \item $ P_i = P(\gamma_i, p_i) $  and $ R_i = R(\gamma_i, p_i), $\item $ f_i = f_{\sigma_i} $ and hence $ R_i =
f_i\circ P_i\circ f_i^{-1}. $\item Let $ \nu_i $ be the Beltrami differentials of
the structures of $ \sigma_i $ respectively.
\endroster
Now our aim is to show that there exists a subsequence $ \{\nu_{i_j}\} $ with norms
uniformly bounded away from 1. By Remark 1 and Lemma 1 it is enough to show that we
can choose a subsequence  $ p_{i_j} $ uniformly bounded away from 1.

Let $ k_i $ be integers so that $ R^{k_i}(\gamma_i) = \gamma. $ Let $ B_i $ be a
component of $ R^{-k_i}(B(\gamma)) $ containing the geodesic $ \gamma_i. $ Then by
our assumptions there exists a constant $ C $ so that the moduli $ m(B_i) =
\frac{m(B(\gamma)}{d(\gamma_i)} \geq C > 0 $ are uniformly bounded. Let $ A_i =
\phi_i(B_i). $

\proclaim{Lemma 2} There exists a subsequence $ \{i_j\} $ and a number $ p < 1 $ so
that $ A(p) \subset A_{i_j} $ for any $j. $
\endproclaim
\demo{Proof} The argument is simple. Let $ A_{i_0} $ be any annulus of  minimal
modulus. Then there exist conformal injections $ h_i :A_{i_0} \mapsto A_i $ such
that $ h_i(\partial\Delta) = \partial\Delta $ and $ h_i(1) = 1. $ The family
$\{h_i\} $ is normal so let $ \{h_{i_j}\} $ be a convergent subsequence. Then the
limit map $ h_\infty \neq const.$ Now let $ q < 1$ be so that $ A(q) \subset
A_{i_0}, $ then by the reflection principle $ \{h_{i_j}\} $ converges to $ h_\infty
$ uniformly on $ A(\frac{q +1}{2}). $ Hence $ h_\infty(A(\frac{q +1}{2})) \subset
h_{i_j}(A_{i_0}) $ for all large enough $ j. $ Let $ p < 1$ be an integer so that $
A(p) \subset h_\infty(A(\frac{q +1}{2})); $ then by the discussion above
$$
A(p) \subset h_{i_j}(A_{i_0}) \subset A_{i_j}
$$
for all large $ j.$ The lemma is thus proved.
\enddemo

By  Remark 1, Lemma 1 and Lemma 2 we have that the family of quasiconformal
homeomorphisms $ \{f_{i_j}\} $ is normal, and after passing to a subsequence we can
assume that $ \{f_{i_j}\} $ converges to a quasi-conformal homeomorphism $
f_\infty.$ The Julia set $ J(R) $ is a Cantor set hence the spherical diameter $
diam(\gamma_i) \to 0. $ Then the homeomorphisms $ f(\gamma_{i_j}, p) $ converge to
the identity uniformly on $ \overline{\Bbb C}, $ and hence $ P_{i_j} \to R. $

Again after passing to a subsequence we can assume that $ lim_{j\to\infty}R_{i_j} =
R_\infty, $ where $ R_\infty $ is a rational map of degree smaller or equal to the
degree of $ R. $ Then we can pass to the limit in the following equality:
$$ f_{i_j}\circ P_{i_j}\circ f^{-1}_{i_j} = R_{i_j} \to f_\infty \circ R\circ
f^{-1}_\infty = R_\infty.
$$
Now to obtain a contradiction assume that $ R $ is a structurally stable map. Then $
R_\infty $ is  structurally stable ( a being a qc-deformation of $ R$) and $ s(R) =
s(R_\infty). $ By  construction  $ R_\infty = lim_{j\to\infty}R_{i_j} $ and by
Corollary 1 $ s(R_{i_j}) \geq s(R) + 1 = s(R_\infty) + 1 $ which contradicts  the
structural stability of $ R_\infty. $

\heading Proof of Theorem 1\endheading

Here we show that the existence of a conical point $ x \in P(R) $ implies Assumption
"G" .
 \proclaim{Lemma 3} Assume $ R $ satisfies  the assumptions of the Theorem 1.
Then $ R $ satisfies  Assumption "G".
\endproclaim
\demo{Proof} Let $ x_0 \in P(R) $ be a conical point. Let integer $ d $ and a
sequence $ \{k_j\} $ be as in the definition of the conical point. Let $ U(x_0,
R^{k_j},\delta) $ be the component of $ R^{-k}({\Bbb D}(R^{k_j}(x_0),\delta)) $ that
contains $x_0,$ where $ {\Bbb D}(R^k(x_0),\delta) $ is the disk centered in $
R^k(x_0) $ with the radius $ \delta.$ After passing to a subsequence we can assume
that the sequence $ R^{k_j}(x_0) $ converges to a point $ y \in P(R).$ The disks $
{\Bbb D}(R^{k_j}(x_0),\delta) $ converges as well as to the disk $ {\Bbb
D}(y,\delta). $ Now let $ \gamma \subset S\cap {\Bbb D}(y,\frac{\delta}{2}) $ with $
y \in I(\gamma). $ Then for large $ j $ we have: \roster\item $ R^{k_j}(x_0) \in
I(\gamma), $ \item there exists a geodesic $ \gamma_j \in R^{-k_j}(\gamma)\cap
U(x_0, R^{k_j},\delta), $ so that $ x_0 \in I(\gamma_j), $ \item the hyperbolic
length $ L(\gamma_j) \leq dL(\gamma). $
\endroster
The lemma is proved.
\enddemo
An application of Theorem A completes the proof of Theorem 1.

\enddocument